\documentclass[12pt,a4,amsmath,amsfonts]{article}
\usepackage{times}
\usepackage{amssymb}
\usepackage{amscd}

\def\Mat{\mathop{\mathrm {Mat}}\nolimits}
\def\Cyl{\mathop{\mathrm {Cyl}}\nolimits}
\def\Cone{\mathop{\mathrm {Cone}}\nolimits}
\def\id{\mathop{\mathrm {id}}\nolimits}
\def\ev{\mathop{\mathrm {ev}}\nolimits}
\newtheorem{thm}{Theorem}[section]
\newtheorem{prop}[thm]{Proposition}
\newtheorem{defn}[thm]{Definition}

\newtheorem{rem}[thm]{Remark}
\begin{document}
\title{Category of Noncommutative CW Complexes\footnote{\bf Version  from \today .
}\footnote{\bf The work was supported in part by Vietnam National
Project for Research in Fundamental Sciences and was completed
during the stay in June - July, 2007 at the Abdus Slam ICTP, Trieste, Italy }}
\author{Do Ngoc Diep}
\maketitle
\begin{abstract}
We expose the notion of noncommutative CW (NCCW) complexes, define noncommutative (NC) mapping cylinder and NC mapping cone, and
prove  the noncommutative Approximation Theorem. The long exact homotopy
sequences associated with  arbitrary morphisms are also
deduced.

\end{abstract}
\section{Introduction}
Classical algebraic topology was fruitfully developed on the
category of topological spaces with CW complex structure, see e.g.
\cite{whitehead}. Our goal is to show that with the same success, theory can be
developed in the framework of noncommutative topology.

In noncommutative geometry the notion of topological
spaces is changed by the notion of C*-algebras, motivating the
spectra of C*-algebras as some noncommutative spaces. In the works \cite{elp} and
\cite{pedersen}, it was introduced the notion of noncommutative CW (NCCW)
complex and proved some elementary properties of NCCW complexes.
We continue this line in proving some basic noncommutative results.
In this paper we aim to explore the same properties of NCCW complexes, as the ones of CW complexes from
algebraic topology. In particular, we prove
some NC Cellular Approximation Theorem and the existence of homotopy exact sequences associated with morphisms. In the work
\cite{diep3} we introduced the notion of NC Serre fibrations
(NCSF) and studied cyclic theories for the (co)homology of these
NCCW complexes. In \cite{diep4} we studied the Leray-Serre spectral
sequences related with cyclic theories: periodic cyclic homology
and KK-theory. In \cite{diepkukutho1} and \cite{diepkukutho2} we
computed some noncommutative Chern characters. Some deep study
should be related with the Busby invariant, studied in
\cite{diep1}, \cite{noncomm}.

Let us describe in more detail the content of the paper. In Section 2 we expose
the pullback and pushout diagrams of G. Pedersen \cite{pedersen}
on categories of C*-algebras. In Section 3 we introduce NCCW
complexes following S. Eilers, T.A. Loring and G. K. Pedersen, etc. We prove in Section 4 a
noncommutative Cellular Approximation Theorem. We prove in Section
5 some long exact homotopy sequences associated with morphisms of
C*-algebras.

\section{Constructions in categories of C*-algebras}
In this section we expose the pullback and pushout constructions of S. Eilers, T.A. Loring and G. K. Pedersen \cite{elp} and  of G. Pedersen
 \cite{pedersen} on categories of  C*-algebras, and after that we define mapping cylinders and mapping cones associatd with arbitrary morphisms.

Let introduce some general notations. By $\mathbf I = [0,1]$
denote the closed interval from $0$ to $1$ on the real line of
real numbers. It is easy to construct a homeomorphism $\mathbf I^n
\approx \mathbf B^n$ between the $n$-cube and the $n$ dimensional
closed ball. Denote also the interior of the cube $\mathbf I^n$ by
$\mathbf I^n_0 = (0,1)^n = \overbrace{\mathbf I^n}^{\circ}$. It is
easy to show that the boundary $\partial \mathbf I^n =\mathbf I^n
\setminus \mathbf I^n_0 $ is homotpic to the $(n-1)$-dimensional
sphere $\mathbf S^{n-1}$. Denote the space of all the continuous
functions on $\mathbf I^n$ with values in a C*-algebra $A$ by
$\mathbf I^n = \mathbf C([0,1]^n,A)$, and by analogy by $\mathbf
I_0^nA := \mathbf C_0((0,1)^n,A)$ the space of all continuous
functions with compact support with values in $A$, and finally, by
$\mathbf S^nA = \mathbf C(\mathbf S^n,A)$ the space of all
continuous maps from $\mathbf S^n$ to $A$.

\begin{defn}[Pullback diagram] {\rm A commutative diagram of C*-algebras and *-homomorphisms
$$\CD X @>\gamma>> B\\
@VV\delta V @VV\beta V\\
A @>\alpha>> C\endCD \eqno{(2.1)}$$
is a} pullback, {\rm if $\ker \gamma \cap \ker \delta = 0$ and if
$$\CD Y @>\psi>> B\\
@VV\varphi V @VV\beta V\\
A @>\alpha>> C\endCD \eqno{(2.2)}$$
is another commutative diagram, then there exists a unique morphism
$\sigma: Y \to X$ such that $\varphi = \delta\circ \sigma$ and
$\psi=\gamma\circ \sigma$, i.e. we have the so called {\it pullback diagram}
\begin{center}
\begin{picture}(100,100)
\put(0,100){Y}
\put(50,50){X}
\put(100,50){B}
\put(50,0){A}
\put(100,0){C}
\put(15,90){\vector(1,-1){30}}
\put(10,90){\vector(1,-2){35}}
\put(15,95){\vector(2,-1){70}}
\put(60,55){\vector(1,0){35}}
\put(60,0){\vector(1,0){35}}
\put(55,45){\vector(0,-1){35}}
\put(105,45){\vector(0,-1){35}}
\put(50,80){$\psi$}
\put(25,80){$\sigma$}
\put(30,30){$\varphi$}
\put(70,58){$\gamma$}
\put(70,10){$\alpha$}
\put(60,30){$\delta$}
\put(110,30){$\beta$}
\put(220,25){(2.3)}
\end{picture}
\end{center}
}\end{defn}

\begin{defn}[Pushout diagram]{\rm A commutative diagram of C*-algebras and *-homomorphisms
$$\CD C @>\beta>> B\\
@VV\alpha V @VV\gamma V\\
A @>\delta>> X\endCD \eqno{(2.4)}$$
is a} pushout, {\rm if $X$ is generated by $ \gamma(B) \cup \delta(A) $  and if
$$\CD C @>>\beta> B\\
@VV\alpha V @VV\psi V\\
A @>\varphi>> Y\endCD \eqno{(2.5)}$$
is another commutative diagram, then there exists a unique morphism $\sigma:
X \to Y$ such that $\varphi = \sigma\circ \gamma$ and $\psi=\sigma\circ \delta$,
 i.e. we have the so called {\it pushout diagram}
\vskip 1cm
\begin{center}
\begin{picture}(100,100)
\put(100,0){Y}
\put(50,50){X}
\put(50,100){B}
\put(0,50){A}
\put(0,100){C}
\put(65,40){\vector(1,-1){30}}
\put(60,90){\vector(1,-2){35}}
\put(15,45){\vector(2,-1){70}}
\put(10,105){\vector(1,0){35}}
\put(10,55){\vector(1,0){35}}
\put(5,95){\vector(0,-1){35}}
\put(55,95){\vector(0,-1){35}}
\put(80,50){$\psi$}
\put(75,30){$\sigma$}
\put(30,30){$\varphi$}
\put(25,58){$\gamma$}
\put(10,70){$\alpha$}
\put(60,70){$\delta$}
\put(25,110){$\beta$}
\put(220,25){(2.6)}
\end{picture}
\end{center}

}
\end{defn}

\begin{defn}[NC cone]{\rm For C*-algebras the} NC cone {\rm of $A$ is defined as the tensor product with } $\mathbf C_0((0,1])$, i.e.
$$\Cone(A) := \mathbf C_0((0,1]) \otimes A.\eqno{(2.7)}$$
\end{defn}

\begin{defn}[NC suspension]{\rm For C*-algebras the} NC suspension {\rm of $A$ is defined as the tensor product with } $\mathbf C_0((0,1))$, i.e.
$$\mathbf S(A) := \mathbf C_0((0,1))\otimes A. \eqno{(2.8)}$$
\end{defn}

\begin{rem}
If $A$ admits a NCCW complex structure, the same have the cone $\Cone(A)$ of $A$ and the suspension $\mathbf S(A)$ of $A$.
\end{rem}

\begin{defn}[NC mapping cylinder]{\rm
Consider a map $f : A \to B$ between C*-algebras. In the algebra $\mathbf C(\mathbf I) \otimes A \oplus B$ consider the closed two-sided ideal $\langle \{1\} \otimes a - f(a), \forall a\in A\rangle $, generated by elements of type $\{1\} \otimes a - f(a), \forall a\in A$. The quotient algebra
$$\Cyl(f)=\Cyl(f:A\to B) := \left(\mathbf C(\mathbf I) \otimes A \oplus B\right)/ \langle \{1\} \otimes a - f(a), \forall a\in A\rangle \eqno{(2.9)}$$ is called the} NC mapping cylinder {\rm and denote it by } $\Cyl(f:A\to B)$.
\end{defn}

\begin{rem}
It is easy to show that $A$ is included in $\Cyl(f: A \to B)$ as $\mathbf C\{0\} \otimes A \subset \Cyl(f:A\to B)$ and $B$ is included in also $B \subset \Cyl(f:A\to B)$.
\end{rem}
\begin{defn}[NC mapping cone]{\rm
In the algebra $\mathbf C((0,1]) \otimes A \oplus B$ consider the closed two-sided ideal $\langle \{1\} \otimes a - f(a), \forall a\in A\rangle $, generated by elements of type $\{1\} \otimes a - f(a), \forall a\in A$.
We define the} mapping cone {\rm as the quotient algebra }
$$\Cone(f)=\Cone(f:A\to B) := \left(\mathbf C_0((0,1]) \otimes A \oplus B\right)/ \langle \{1\} \otimes a - f(a), \forall a\in A\rangle .\eqno{(2.10)}$$
\end{defn}
\begin{rem}
It is easy to show that $B$ is included in $\Cone(f: A \to B)$.
\end{rem}

\begin{prop}\label{diags}
Both the mapping cylinder and mapping cone satisfy the pullback
diagrams
$$\CD \Cone(\varphi) @>pr_1 >> \mathbf C_0(0,1] \otimes A \\
@V{pr_2}VV  @VV{\varphi\circ \ev(1)}V\\
B @>id>> B\endCD \qquad
\CD \Cyl(\varphi) @>pr_1 >> \mathbf C[0,1] \otimes A \\
@V{pr_2}VV  @VV{\varphi\circ \ev(1)}V\\
B @>id>> B\endCD,$$ where $ev(1)$ is the map of evaluation at the
point $1\in [0,1]$.
\end{prop}
\begin{rem}
The pullback diagrams in Proposition \ref{diags} can be used as the initial definition of mapping cyliner and mapping cone. The previous definitions are therefore the existence of those universal objects.
\end{rem}

\begin{rem} It is reasonable to have that
in the case of C*-algebras of continuous functions $A=\mathbf C(X)$,
$B=\mathbf C(Y)$ the C*-algebras of continuous functions over the
cone and the suspension of topological spaces are in general
different from the cone and the suspension of C*-algebras, we have
just defined; the same is true that the mapping cylinder and the
mapping cone of morphisms of C*-algebras are different from the
C*-algebra of continuous functions on the mapping cylinder and the
mapping cone of spaces.
\end{rem}

\section{The category NCCW}
In this section we introduce NCCW complexes following J. Cuntz and
following S. Eilers, T. A. Loring and G. K. Pedersen, \cite{elp}
etc.
\begin{defn} A dimension 0 NCCW complex {\rm is defined, following
\cite{pedersen} as a finite sum of C* algebras of finite linear
dimension, i.e. a sum of finite dimensional matrix algebras,}
$$A_0 = \bigoplus_{k} \mathbf M_{n(k)}.\eqno{(3.1)}$$
In dimension n, an NCCW complex  {\rm is defined
 as a sequence $\{A_0,A_1,\dots,A_n\}$ of C*-algebras $A_k$
 obtained each from the previous one by the pullback construction
 $$\CD 0 @>>> I^k_0F_k @>>> A_k @>\pi>> A_{k-1} @>>> 0\\
 @. @| @VV\rho_k V  @VV\sigma_k V @. \\
 0 @>>> I_0^kF_k @>>> I^kF_k @>\partial>> \mathbf S^{k-1}F_k @>>>
 0,\endCD \eqno{(3.2)}$$ where $F_k$ is some C*-algebra of finite
 linear dimension, $\partial$ the restriction morphism, $\sigma_k$
 the connecting morphism, $\rho_k$ the projection on the first
 coordinates and $\pi$ the projection on the second coordinates in
 the presentation
 $$A_k = \mathbf I^kF_k \bigoplus_{\mathbf S^{k-1}F_k}A_{k-1}\eqno{(3.3)}$$
 }\end{defn}
\begin{prop} \label{cyl}
If the algebras $A$ and $B$ admit a NCCW complex structure, then
the same has the NC mapping cylinder $\Cyl(f:A\to B)$.
\end{prop}
{\sc Proof}.  Let us remember from \cite{pedersen} that the interval $\mathbf I$ admits a structure of an NCCW complex. Next, tensor product of two NCCW complex \cite{pedersen} is also an NCCW complex and finally the quotient of an NCCW complex is also an NCCW complex, loc. cit.. \hfill$\Box$
\begin{prop}
If the algebras $A$ and $B$ admit a NCCW complex structure, then the same has the NC mapping  cone $\Cone(f:A\to B)$.
\end{prop}
{\sc Proof}. The same argument as in Proof of Proposition \ref{cyl}. \hfill$\Box$

\section{Approximation Theorem}
We prove in this section a noncommutative analog of the well-known
Cellular Approximation Theorem. First we introduce the so called
noncommutative homotopy extension property (NC HEP).
\begin{defn}[NC HEP] For a given $(f,\varphi_t)$ and a C*-algebra $C$, we say that $\tilde{h}=\tilde{\varphi}_t$ is
a solution of the
extension problem if we have the commutative homotopy extension
diagram
\begin{center}
\begin{picture}(200,50)
\put(195,50){C} \put(0,0){$\Cyl(i: A \hookrightarrow B)$}
\put(160,0){$C[0,1] \otimes B$} \put(200,45){\vector(0,-1){30}}
\put(190,45){\vector(-4,-1){130}}
\put(140,5){\vector(-1,0){50}} \put(50,25){$(f,\varphi_t)$}
\put(185,25){$\tilde{\varphi_t}$}
\end{picture}
\end{center}
\end{defn}

\begin{defn}[NC NDR] We say that the pair of algebras $(B,A)$ is a
NCNDR pair, if there are continuous morphisms $u: \mathbf
C[0,1]\to  B $ and $\varphi :  B\to\mathbf C[0,1] \otimes B\cong \mathbf
C(I,B) $ such that
\begin{enumerate}
\item $u^{-1}(A)=0$;
\item If $\varphi(b) = (x(t),b')$ and $x(t)=0\in \mathbf C(\mathbf I)$ then $b'=b$, $\forall b\in B$;
\item $  \varphi(a)=(x(t),a), \forall a\in A,  x(t)\in\mathbf
C(\mathbf I)$;
\item $\varphi(b) = (x(t),b')$ and if $x(t)=1\in \mathbf C(\mathbf I)$ then $b'\in A$ for all $b\in B$ such that $u(b) \neq 1$.
\end{enumerate}
\end{defn}
The following proposition is easily to prove.
\begin{prop}
The assertion that NC HEP has solution for every $\varphi_t$ and
$C$ is equivalent to the the property that $(B,A)$ is a NC NDR
pair.
\end{prop}
{\sc Proof}. If $(B,A)$ has NC HEP, we can for every $C$,
construct $\tilde{\varphi} :  B\to\mathbf C[0,1] \otimes B $
satisfying the NC HEP diagram. Choose $C=B$ and $f=\id$ we have
the function $\varphi$ and then choose $(D,C) = (\mathbf C[0,1],0)$
in the definition of NC NDR pair we have the function $u$.

Conversely, if $(B,A)$ is a NC NDR pair, we can define
$$h=\varphi : B\to\mathbf C[0,1] \otimes B  ,\eqno{(4.1)}$$
the composition of which with $f:C\to B  $ satisfy the NC HEP
diagram. \hfill $\Box$

\begin{thm}[Extension]\label{ext} Suppose that $B = \mathbf I^nF_n \oplus_{\mathbf S^{n-1}F_n} A$
and $(C,D)$ is a NC NDR pair. Every relative morphism of pairs of
C*-algebras
$$f:(D,C)\to(\Cyl(i:A\hookrightarrow B), \mathbf C\{1\} \otimes A)  \eqno{(4.2)}$$ can be up-to homotopy extended to a relative morphism of pairs of C*-algebras
$$F: (D,C)\to(\mathbf C(\mathbf I) \otimes B, \mathbf C\{1\} \otimes B)  .\eqno{(4.3)}$$
\end{thm}
{\sc Proof}. The property that $(D,C)$ is a NC NDR pair, there is
a natural extension
$$f_1:(D,C)\to(\mathbf C(\mathbf I)\otimes B,\mathbf C\{1\} \otimes A)  .\eqno{(4.4)}$$
Composing $f_1$ with the map, evaluating the value at $1$  give a
morphism
$$\ev(1)\circ f_1 : (D,C)\to(\mathbf C\{1\} \otimes B,\mathbf C\{1\} \otimes A)  .\eqno{(4.5)}$$
Therefore, there exists a natural extension $f_2$ from the pair $(D,C).$ to the pair
$$(\mathbf C\{0\} \otimes \mathbf C(\mathbf I) \otimes B + \mathbf C(\mathbf I) \otimes \mathbf C\{1\} \otimes B, \Cyl(\mathbf C\{1\} \otimes A \hookrightarrow \mathbf C\{1\} \otimes B)).$$
Once again, there is a natural extension $f_3$ from the pair to the pair
$$(\mathbf C\{0\} \otimes \mathbf C(\mathbf I) \otimes B + \mathbf C\{1\} \otimes \mathbf C(\mathbf I) \otimes B + \mathbf C(\mathbf I) \otimes \mathbf C(\mathbf I) \otimes A, Cyl(\mathbf C(\mathbf I) \otimes A \hookrightarrow \mathbf C(\mathbf I) \otimes B))=$$
$$=(\mathbf C\{0\} \otimes \Cyl(\mathbf C(\mathbf I) \otimes A \hookrightarrow \mathbf C(\mathbf I) \otimes B) + \mathbf C\{1\} \otimes \Cyl(\mathbf C(\mathbf I) \otimes A \hookrightarrow \mathbf C(\mathbf I) \otimes B), $$ $$,\Cyl(\mathbf C(\mathbf I) \otimes A \hookrightarrow \mathbf C(\mathbf I) \otimes B)).$$

And finally, there is a natural extension $f_4$ from the pair $(D,C)$ to the pair
$$\mathbf C(\mathbf I) \otimes \mathbf C(\mathbf I) \otimes B, \Cyl(\mathbf C(\mathbf I) \otimes A \hookrightarrow \mathbf C(\mathbf I) \otimes B)). $$ We define the desired extension $$F:  (D,C)\to(\mathbf C(\mathbf I) \otimes B, \mathbf C\{1\} \otimes B) $$ as
$$F(t,x) := f_4(t,0,x).\eqno{(4.6)}$$
\hfill$\Box$

\begin{thm}
Let $\{ A_0,A_1,\dots,A_n\}$  and $\{B_0, B_1,\dots,B_m\}$ be two
NCCW complexes and $f: A = A_n \to B_m = B$ an algebraic
homomorphism (map). Then $f$ is homotopic to a cellular NCCW
complex  map  $h: A \to B$.
\end{thm}
{\sc Proof}. We constract a sequence of maps
$$ g_p :\mathbf  A_p \to C(\mathbf I) \otimes B_p, \eqno{(4.7)}$$
with 4 well-known properties:
\begin{enumerate}
\item $g_p(x) = (0,f(x)), \forall x\in A_p$
\item If $g_p(b) = (x(t),f(b))$ and if $x(t) = 0\in\mathbf C(\mathbf I)$, then
$f(b) = b, \forall b\in B$.
\item $\ev(1)\circ g_p = g_{p-1},$
\item $g_p( A_p) \subset \mathbf C\{1\} \otimes B_p .$
\end{enumerate}
Indeed, following the definition of an NCCW complex structure, we
have
$$A_0 = F_0 \otimes A,\qquad
F_0 = \bigoplus_{j_0} \mathbf M_n(j_0)\eqno{(4.8)}$$ a finite
system of quantum points,  i.e. a commutative diagram
$$\CD 0 @>>>\mathbf I_0^1 F_1 @>>> A_1 @>>> A_0 @>>> 0\\
@. @| @VV\rho_1 V  @VV\sigma_1 V @.\\
0 @>>> \mathbf I_0^1F_1 @>>> \mathbf I^1F_1 @>>\alpha_1> \mathbf
S^0F_1 @>>> 0\endCD \eqno{(4.9)}$$ in which the second square is a
pullback diagram,
$$F_1 = \bigoplus_{j_1} \mathbf M_{n(j_1)} =
\bigoplus_{j_1}\Mat_{n(j_1)},\eqno{(4.10)}$$ and we can present
$A_1$ as
$$ A_1 \approx \mathbf I^1F_1 \bigoplus_{S^0F_1} A_0. \eqno{(4.11)}$$
Following the compessible theorems \cite{whitehead} and the previous Extension Theorem \ref{ext}, the function
$g_0$ can be naturally extended to a function $g_1$ with
properties 1. - 4. and now we have again following the definition
of an NCCW complex,
$$\CD 0 @>>>\mathbf I_0^2 F_2 @>>> A_2 @>>> A_1 @>>> 0\\
@. @| @VV\rho_2 V  @VV\sigma_2 V @.\\
0 @>>> \mathbf I_0^2F_2 @>>> \mathbf I^2F_2 @>>\alpha_2> \mathbf
S^1F_2 @>>> 0\endCD \eqno{(4.12)}$$
$$F_2 = \bigoplus_{j_2} \mathbf M_{n(j_2)} =
\bigoplus_{j_2}\Mat_{n(j_2)},\eqno{(4.13)}$$ and we can present
$A_2$ as
$$ A_2 \approx \mathbf I^2F_2 \bigoplus_{S^1F_2} A_1. \eqno{(4.14)}$$
Following the compessible theorems \cite{whitehead} and the previous Extension Theorem \ref{ext}, the function
$g_1$ can be naturally extended to a function $g_2$ with
properties 1. - 4. The procedure is continued for all $p$. Once
these functions $g_p$ were defined, the function $g:  A \to \mathbf
C(\mathbf I) \otimes B$ which is continuous and $g$ is a
homotopy of $f$ to $h$, where $$h(x) := \ev(1)\circ g(x).\eqno{(4.15)}$$
Because of 4. the function $h : A \to B$ is a cellular NCCW
complex map. \hfill$\Box$

\section{Homotopy of NCCW complexes}
We prove in this section the standard long exact homotopy
sequences.

Let us first recall the definition of homotopic morphisms.
\begin{defn}
A homotopy {\rm between two morphisms $\varphi, \psi : A \to B$ is
a morphism $\Phi : A \to \mathbf C(\mathbf I) \otimes B$, such
that $\Phi(0,.) = \varphi$ and $\Phi(1,.) = \psi$. }
\end{defn}

\begin{prop}\label{homot}
There is a natural homotopy $\Cyl(\varphi : A \to B) \simeq B$ and
$\Cone(\varphi : A \to B)\simeq B/A$, if the last one $B/A$ is
defined.
\end{prop}

\begin{thm}
For every morphism $\varphi ; A \to B$, there is a natural long
exact homotopy sequence
$$\CD \dots @>>> \mathbf S^2(A) @>>> \mathbf S(\Cone(\varphi : A \to B)) @>>>\mathbf
S(\Cyl(\varphi : A \to B)) @>>>\endCD$$
$$\CD S(A) @>>> \Cone(\varphi : A \to B) @>>> \Cyl(\varphi : A \to B) @>>> A @>\varphi>> B\endCD\eqno{(5.1)}$$
\end{thm}
{\sc Proof}. Put $A_0 = B, A_1= A$ and $\varphi_0 = \varphi$ we
have
$$\CD A_0 @>\varphi_0 = \varphi>> A_1 .\endCD\eqno{(5.2)}$$
Because of Proposition \ref{homot} we have
$$\CD A_0=B @<\varphi_0<< A_1= A @<\varphi_1<< A_2=
\Cyl(\varphi) @<\varphi_2<< A_3=\Cone(\varphi).\endCD
\eqno{(5.3)}$$
Because of the exact sequence
$$\CD\mathbf S(A)@<<<\Cyl(\varphi)  @<<< \Cone(\varphi)  ,\endCD$$
we have
$$\CD A_0=B @<\varphi_0<< A_1= A @<\varphi_1<< A_2=
\Cyl(\varphi) @<\varphi_2<<
A_3=\Cone(\varphi)@<\varphi_3<<\endCD$$
$$\CD @<\varphi_3<< A_4= \mathbf S(A)= \mathbf C_0((0,1)) \otimes A. \endCD\eqno{(5.4)}$$
Because the tensor product $\mathbf C_0((0,1)) \otimes .$ is a
left exact functor and because of (5.3) , we have
$$\CD A_0=B @<\varphi= \varphi_0<< A_1=A @<\varphi_1<< A_2 = \Cyl(\varphi)
@<\varphi_2<< A_3= \Cone(\varphi) @<\varphi_3<<\endCD$$
$$\CD@<\varphi_3<< \mathbf S(A) @<\varphi_4<<
  \mathbf S(\Cyl(\varphi)) @<\varphi_5 <<
A_5=\mathbf S(\Cone(\varphi)) @<\varphi_5 << \endCD$$
$$\CD @<\varphi_5 <<A_6=\mathbf S^2(A)@<\varphi_6<<\dots ,
\endCD\eqno{(5.5)}$$ etc. \hfill$\Box$

\section*{Acknowledgments}
The work was supported in part by Vietnam National Project for
Research in Fundamental Sciences and was completed  during the
stay in June and July, 2007 of the author, in Abdus Salam ICTP,
Trieste, Italy.  The author expresses his deep and sincere thanks
to Abdus Salam ICTP and especially Professor Dr. Le Dung Trang for
the invitation and for providing the nice conditions of work, and
Professor C. Schochet for some Email discussions.

\vskip 0.5cm
\parindent=0pt
{\sc Institute of Mathematics, Vietnam Academy of Science and Technology, 18 Hoang Quoc Viet Road, Cau Giay District, 10307 Hanoi, Vietnam}\\
{\tt Email: dndiep@math.ac.vn}

\end{document}